# WHEN THE LAW OF LARGE NUMBERS FAILS FOR INCREASING SUBSEQUENCES OF RANDOM PERMUTATIONS

By Ross G. Pinsky

*Technion—Israel Institute of Technology*

Let the random variable $Z_{n,k}$ denote the number of increasing subsequences of length $k$ in a random permutation from $S_n$, the symmetric group of permutations of $\{1,\ldots,n\}$. In a recent paper [*Random Structures Algorithms* **29** (2006) 277–295] we showed that the weak law of large numbers holds for $Z_{n,k_n}$ if $k_n = o(n^{2/5})$; that is,

$$\lim_{n\to\infty} \frac{Z_{n,k_n}}{EZ_{n,k_n}} = 1 \qquad \text{in probability.}$$

The method of proof employed there used the second moment method and demonstrated that this method cannot work if the condition $k_n = o(n^{2/5})$ does not hold. It follows from results concerning the longest increasing subsequence of a random permutation that the law of large numbers cannot hold for $Z_{n,k_n}$ if $k_n \geq cn^{1/2}$, with $c > 2$. Presumably there is a critical exponent $l_0$ such that the law of large numbers holds if $k_n = O(n^l)$, with $l < l_0$, and does not hold if $\limsup_{n\to\infty} \frac{k_n}{n^l} > 0$, for some $l > l_0$. Several phase transitions concerning increasing subsequences occur at $l = \frac{1}{2}$, and these would suggest that $l_0 = \frac{1}{2}$. However, in this paper, we show that the law of large numbers fails for $Z_{n,k_n}$ if $\limsup_{n\to\infty} \frac{k_n}{n^{4/9}} = \infty$. Thus, the critical exponent, if it exists, must satisfy $l_0 \in [\frac{2}{5}, \frac{4}{9}]$.

**1. Introduction and statement of results.** Let $S_n$ denote the symmetric group of permutations of $\{1,\ldots,n\}$. By introducing the uniform probability measure $U_n$ on $S_n$, one can consider $\sigma \in S_n$ as a random permutation. Probabilities and expectations according to $U_n$ will frequently be denoted by the notation $P$ and $E$, respectively. The problem of analyzing the distribution of the length, $L_n$, of the longest increasing subsequence in a random permutation from $S_n$ has a long and distinguished history; see [1] and references therein. In particular, the work of Logan and Shepp [3] together with that of Vershik and Kerov [5] show that $EL_n \sim 2n^{1/2}$ and that $\sigma^2(L_n) = o(n)$, as $n \to \infty$. More recently, the profound work of Baik, Deift and Johansson [2] has led to the following theorem.

THEOREM BDJ.
$$\lim_{n\to\infty} P\left(\frac{L_n - 2n^{1/2}}{n^{1/6}} \leq x\right) = F(x),$$

*where $F$ is the Tracy–Widom distribution.*

Consider now the random variable $Z_{n,k} = Z_{n,k}(\sigma)$, which we define to be the number of increasing subsequences of length $k$ in a permutation $\sigma \in S_n$. Thus, for example, if $\sigma = \begin{pmatrix} 1\ 2\ 3\ 4\ 5 \\ 1\ 3\ 4\ 5\ 2 \end{pmatrix}$, then $Z_{5,3}(\sigma) = 4$ because there are four increasing subsequences of length three; namely, 134, 135, 145 and 345. It is useful to represent $Z_{n,k}$ as a sum of indicator random variables. For positive integers $\{x_1, \ldots, x_k\}$ satisfying $1 \leq x_1 < x_2 < \cdots < x_k \leq n$, let $B^n_{x_1,\ldots,x_k} \subset S_n$ denote the subset of permutations which contain the increasing subsequence $\{x_1, x_2, \ldots, x_k\}$. Then we have

(1.1) $$Z_{n,k} = \sum_{x_1 < x_2 < \cdots < x_k} 1_{B^n_{x_1, x_2, \ldots, x_k}},$$

where the sum is over the $\binom{n}{k}$ distinct increasing subsequences of length $k$. Since the probability that a random permutation fixes any particular increasing sequence of length $k$ is $\frac{1}{k!}$, it follows that the expected value of $Z_{n,k}$ is given by

(1.2) $$EZ_{n,k} = \frac{\binom{n}{k}}{k!}.$$

One can consider $k$ to depend on $n$ in which case we write $k_n$. A straightforward calculation using Stirling's formula shows that

(1.3)
$$EZ_{n,cn^l} \sim \frac{1}{2\pi cn^l}\left[\left(\frac{e}{c}\right)^2 n^{1-2l}\right]^{cn^l}, \qquad \text{as } n \to \infty, \text{ for } l \in \left(0, \frac{1}{2}\right);$$
$$EZ_{n,cn^{1/2}} \sim \frac{\exp(-c^2/2)}{2\pi cn^{1/2}}\left(\frac{e}{c}\right)^{2cn^{1/2}}, \qquad \text{as } n \to \infty.$$

In [4] we proved the following law of large numbers for $Z_{n,k_n}$.

THEOREM P. *Let $k_n = o(n^{2/5})$. Then*

(1.4) $$\lim_{n\to\infty} \frac{Z_{n,k_n}}{EZ_{n,k_n}} = 1 \qquad \text{in probability.}$$

We note that the proof in [4] used the second moment method. In fact, it was shown that $\text{Var}(Z_{n,k_n}) = o((EZ_{n,k_n})^2)$ if and only if $k_n = o(n^{2/5})$. Although we do not have a proof of this, presumably there exists a critical exponent $l_0$ such that the law of large numbers holds for $Z_{n,n^l}$ when $l < l_0$ and does not hold when $l > l_0$. Note that several phase transitions occur when $l = \frac{1}{2}$. First, from [3] and [5], or from [2], it follows that $\lim_{n\to\infty} P(Z_{n;cn^{1/2}} = 0) = 1$ or 0, depending on whether $c > 2$ or $c < 2$; in particular then, it is not possible for the law of large numbers to hold for $Z_{n,cn^{1/2}}$ with $c > 2$. Consequently, if a critical exponent exists, it must be less than or equal to $\frac{1}{2}$. Second, from (1.3), it follows that $\lim_{n\to\infty} EZ_{n,cn^{1/2}} = \infty$, if $c < e$, and $\lim_{n\to\infty} EZ_{n,cn^{1/2}} = 0$, if $c \geq e$. And third, note that the factor $\exp(-\frac{c^2}{2})$ suddenly appears in the formula for $EZ_{n,cn^l}$ in (1.3) when $l = \frac{1}{2}$. In light of these phase transitions, we felt that the critical exponent was probably equal to $\frac{1}{2}$. Thus, the main result of this paper came to us as a surprise.

THEOREM 1. *The law of large numbers does not hold for $Z_{n,k_n}$ if*

$$\limsup_{n\to\infty} \frac{k_n}{n^{4/9}} = \infty.$$

Thus, assuming that the critical exponent exists, its value has now been narrowed down to the closed interval $[\frac{2}{5}, \frac{4}{9}]$.

In [4] it was shown that the law of large numbers in (1.4) is equivalent to a certain approximation result for the uniform measure on $S_n$, which we now describe. Recall that for probability measures $P_1$ and $P_2$ on $S_n$, the total variation norm is defined by

$$\|P_1 - P_2\| \equiv \max_{A \subset S_n} (P_1(A) - P_2(A)) = \frac{1}{2} \sum_{\sigma \in S_n} |P_1(\sigma) - P_2(\sigma)|.$$

For $x_1 < x_2 < \cdots < x_{k_n}$, let $U_{n;x_1,x_2,\ldots,x_{k_n}}$ denote the uniform measure on permutations which have $\{x_1, x_2, \ldots, x_{k_n}\}$ as an increasing sequence; that is, $U_{n;x_1,x_2,\ldots,x_{k_n}}$ is uniform on $B^n_{x_1,x_2,\ldots,x_{k_n}}$. Note that $U_{n;x_1,x_2,\ldots,x_{k_n}}$ is defined by $U_{n;x_1,x_2,\ldots,x_{k_n}}(\sigma) = \frac{k_n!}{n!} 1_{B^n_{x_1,x_2,\ldots,x_{k_n}}}(\sigma)$. Now define the probability measure $\mu_{n;k_n}$ on $S_n$ by

(1.5) $$\mu_{n;k_n} = \frac{1}{\binom{n}{k_n}} \sum_{x_1 < x_2 < \cdots < x_{k_n}} U_{n;x_1,x_2,\ldots,x_{k_n}}.$$

PROPOSITION P. *The law of large numbers holds for $Z_{n,k_n}$; that is,*

$$\lim_{n\to\infty} \frac{Z_{n,k_n}}{EZ_{n,k_n}} = 1 \quad \text{in probability,}$$

*if and only if*

(1.6) $$\lim_{n\to\infty} \|\mu_{n;k_n} - U_n\| = 0.$$

REMARK. The measure $\mu_{n;k_n}$ can be realized in the following way. Lay down in random order a row of $n$ cards numbered from 1 to n. Choose $k_n$ cards at random from the row, remove them and then replace them in the $k_n$ vacant places in increasing order. The resulting distribution is $\mu_{n;k_n}$. In light of this realization of $\mu_{n;k_n}$, Theorem P and Theorem 1 along with Proposition P give the following interesting interpretation. Take a random deck of cards and "adulterate" it by removing at random $k_n$ cards and then replacing them in increasing order. Then asymptotically, as $n \to \infty$, it is possible to detect this adulteration if $k_n$ is of an order larger than $n^{4/9}$, while it is not possible to detect the adulteration if $k_n = o(n^{2/5})$.

The proof of Theorem 1 exploits Proposition P and Theorem BDJ.

**2. Proof of Theorem 1.** Recall that $L_n = L_n(\sigma)$ denotes the length of the longest increasing subsequence in the permutation $\sigma \in S_n$. Theorem 1 follows easily from the following lemma along with Theorem BDJ and Proposition P.

LEMMA 1. *If* $\limsup_{n\to\infty} \frac{k_n}{n^{1/2}} < \infty$ *and* $\limsup_{n\to\infty} \frac{k_n}{n^{4/9}} = \infty$, *then there exists a sequence* $\{c_n\}_{n=1}^{\infty}$ *satisfying* $\lim_{n\to\infty} c_n = \infty$ *and such that*

$$\limsup_{n\to\infty} \mu_{n,k_n}(L_n > 2n^{1/2} + c_n n^{1/6}) > 0. \tag{2.1}$$

PROOF OF THEOREM 1. Assume to the contrary that the law of large numbers holds for some $\{k_n\}_{n=1}^{\infty}$ satisfying $\limsup_{n\to\infty} \frac{k_n}{n^{4/9}} = \infty$. By Theorem BDJ (or even by the weaker results of Logan and Shepp [3] and of Vershik and Kerov [5]), we may assume that $\limsup_{n\to\infty} \frac{k_n}{n^{1/2}} < \infty$. It then follows from Proposition P and Lemma 1 that $\limsup_{n\to\infty} P(L_n > 2n^{1/2} + c_n n^{1/6}) > 0$. But this contradicts Theorem BDJ. □

It remains to prove Lemma 1.

PROOF OF LEMMA 1. By considering a subsequence if necessary, we may assume without loss of generality that $\lim_{n\to\infty} \frac{k_n}{n^{4/9}} = \infty$. Under this assumption, we will show that (2.1) holds with lim sup replaced by lim inf. Consider a permutation $\sigma \in S_n$ chosen according to the measure $\mu_{n;k_n}$, defined above in (1.5). The permutation $\sigma$ is obtained in two steps as follows. First choose randomly (i.e., according to the uniform distribution) one of the $\binom{n}{k_n}$ increasing subsequences $\{x_1, x_2, \ldots, x_{k_n}\}$ from $\{1, 2, \ldots, n\}$. Then choose randomly one of the permutations that have $\{x_1, x_2, \ldots, x_{k_n}\}$ as an increasing sequence; that is, choose one of the permutations from $B^n_{x_1,x_2,\ldots,x_{k_n}}$ according to the uniform measure $U_{n;x_1,x_2,\ldots,x_{k_n}}$. Note that as a consequence of the second step, the collection of $k_n$ positions occupied by $\{x_1, x_2, \ldots, x_{k_n}\}$ is random.

Let $r_n = n - k_n$ and, having fixed an increasing subsequence $\{x_1, x_2, \ldots, x_{k_n}\}$, let $\{y_1, y_2, \ldots, y_{r_n}\}$ denote the increasing subsequence consisting of all the elements of $\{1, 2, \ldots, n\} - \{x_1, x_2, \ldots, x_{k_n}\}$. Let $L_{n;y_1,y_2,\ldots,y_{r_n}} = L_{n;y_1,y_2,\ldots,y_{r_n}}(\sigma)$ denote the length of the longest increasing subsequence in $\sigma \in B^n_{x_1,x_2,\ldots,x_{k_n}}$ that can be constructed from the elements $\{y_1, y_2, \ldots, y_{r_n}\}$. Fix a sequence $\{\gamma_n\}_{n=1}^\infty$ satisfying $\lim_{n\to\infty} \gamma_n = \infty$. Since $U_{n;x_1,x_2,\ldots,x_{k_n}}$ is uniform on $B^n_{x_1,x_2,\ldots,x_{k_n}}$, it follows from Theorem BDJ that

$$(2.2) \quad \lim_{n\to\infty} U_{n;x_1,x_2,\ldots,x_{k_n}}(L_{n;y_1,y_2,\ldots,y_{r_n}} \geq 2r_n^{1/2} - \gamma_n r_n^{1/6}) = 1.$$

Since $\limsup_{n\to\infty} \frac{k_n}{n^{1/2}} < \infty$, we have

$$(2.3) \quad \begin{aligned} 2r_n^{1/2} - \gamma_n r_n^{1/6} &= 2(n-k_n)^{1/2} - \gamma_n(n-k_n)^{1/6} \\ &\geq 2n^{1/2} - 2\gamma_n n^{1/6} \quad \text{for large } n. \end{aligned}$$

From (2.2) and (2.3) we conclude that

$$(2.4) \quad \lim_{n\to\infty} U_{n;x_1,x_2,\ldots,x_{k_n}}(L_{n;y_1,y_2,\ldots,y_{r_n}} \geq 2n^{1/2} - 2\gamma_n n^{1/6}) = 1.$$

Now define

$$(2.5) \quad \begin{aligned} A_{n;y_1,y_2,\ldots,y_{r_n}} = \{\sigma \in B^n_{x_1,x_2,\ldots,x_{k_n}} : &\sigma \text{ possesses an increasing subsequence} \\ &\text{of length } [2n^{1/2} - 2\gamma_n n^{1/6}] \text{ constructed} \\ &\text{from elements of } \{y_1, y_2, \ldots, y_{r_n}\}\}. \end{aligned}$$

By (2.4), we have

$$(2.6) \quad \lim_{n\to\infty} U_{n;x_1,x_2,\ldots,x_{k_n}}(A_{n;y_1,y_2,\ldots,y_{r_n}}) = 1.$$

Let

$$V_{n;x_1,x_2,\ldots,x_{k_n}}(\cdot) = U_{n;x_1,x_2,\ldots,x_{k_n}}(\cdot \mid A_{n;y_1,y_2,\ldots,y_{r_n}}).$$

Of course, $V_{n;x_1,x_2,\ldots,x_{k_n}}$ is uniform on $A_{n;y_1,y_2,\ldots,y_{r_n}}$. Let

$$s_n = [2n^{1/2} - 2\gamma_n n^{1/6}].$$

Note that from the uniformity it follows that the $V_{n;x_1,x_2,\ldots,x_{k_n}}$-probability that a particular increasing subsequence $\{z_1, z_2, \ldots, z_{s_n}\}$ from $\{y_1, y_2, \ldots, y_{r_n}\}$ belongs to $\sigma \in A_{n;y_1,y_2,\ldots,y_{r_n}}$ and appears in a particular collection of $s_n$ positions is the same for every such subsequence and every such collection of $s_n$ positions.

We summarize the above analysis as follows. Under $\mu_{n;k_n}$, a permutation $\sigma \in S_n$ will have an increasing subsequence $\{x_1, x_2, \ldots, x_{k_n}\}$ of length $k_n$, chosen randomly from all such increasing subsequences, and positioned randomly. Furthermore, in light of (2.6), with high probability it will also have an increasing subsequence $\{z_1, z_2, \ldots, z_{s_n}\}$ of length $s_n$ and disjoint from $\{x_1, x_2, \ldots, x_{k_n}\}$. Given that such a second increasing subsequence exists, the probability that any particular such increasing subsequence exists and occupies any particular collection of $s_n$ positions is the same for all such particular increasing subsequences and all such particular collections of $s_n$ positions.

We now wish to determine how long an increasing subsequence can be constructed from the two increasing subsequences $\{x_1, x_2, \ldots, x_{k_n}\}$ and $\{z_1, z_2, \ldots, z_{s_n}\}$. In light of the above analysis, we consider a row of $n$ spaces and a collection of $n$ cards numbered from 1 to $n$. We randomly select $k_n$ cards and randomly determine $k_n$ spaces in which to place these cards in increasing order from left to right. Then we randomly select $s_n$ other cards and $s_n$ other spaces in which to place these other cards in increasing order from left to right. Note that the positions of the set of $k_n$ cards and the set of $s_n$ cards relative to each other are random and thus independent of $n$. With regard to the numbers appearing on the set of $k_n$ cards and the set of $s_n$ cards, note that what is relevant for determining increasing subsequences is the relative rankings of these numbers rather than the numbers themselves. In light of these facts, if follows that the number $n$ serves no purpose. We might as well consider the above setup with exactly $k_n + s_n$ cards and spaces, numbered from 1 to $k_n + s_n$. We will prove the following lemma.

LEMMA 2. *Let $s_n$ satisfy*

$$c_1 n^{1/2} \leq s_n \leq c_2 n^{1/2},$$

*where $c_1, c_2 > 0$, and let $k_n$ satisfy*

$$\lim_{n \to \infty} \frac{k_n}{n^{1/3}} = \infty \quad \text{and} \quad \limsup_{n \to \infty} \frac{k_n}{n^{1/2}} < \infty.$$

*Consider a row of $s_n + k_n$ spaces and consider a collection of $s_n + k_n$ cards numbered from 1 to $s_n + k_n$. Randomly select $k_n$ cards and randomly select $k_n$ spaces. Insert the $k_n$ cards into the $k_n$ spaces in increasing order from left to right. Now insert the remaining $s_n$ cards in the remaining $s_n$ spaces in increasing order from left to right. Then for some $\delta > 0$,*

$$\liminf_{n \to \infty} \text{Prob}\left(\text{there exists an increasing subsequence of length } s_n + \left[\delta \frac{k_n^{3/2}}{n^{1/2}}\right]\right) > 0.$$

We now complete the proof of Lemma 1. Let $k_n$ be as in the statement of Lemma 1. Then $k_n$ also satisfies the requirement of Lemma 2. Recall that (2.5) and (2.6) hold for any choice of $\{\gamma_n\}$ satisfying $\lim_{n \to \infty} \gamma_n = \infty$.

By the assumption in Lemma 1, $\lim_{n\to\infty} \frac{k_n}{n^{4/9}} = \infty$; thus, $\lim_{n\to\infty} \frac{k_n^{3/2}}{n^{2/3}} = \infty$. Therefore, $[\delta \frac{k_n^{3/2}}{n^{1/2}}] = d_n n^{1/6}$, where $\lim_{n\to\infty} d_n = \infty$. We apply Lemma 2 with $k_n$ as above and with $s_n = [2n^{1/2} - 2\gamma_n n^{1/6}]$, where we choose $\gamma_n$ such that $\lim_{n\to\infty} \gamma_n = \infty$ and $\gamma_n \leq \frac{d_n}{4}$. It follows from the above analysis of $\mu_{n;k_n}$ and from (2.6) that

$$(2.7) \qquad \liminf_{n\to\infty} \mu_{n;k_n}(L_n \geq [2n^{1/2} - 2\gamma_n n^{1/6}] + d_n n^{1/6}) > 0.$$

Define $c_n = \frac{d_n}{3}$. Since $[2n^{1/2} - 2\gamma_n n^{1/6}] + d_n n^{1/6} \geq 2n^{1/2} + \frac{d_n}{2} n^{1/6} - 1 \geq 2n^{1/2} + c_n n^{1/6}$, Lemma 1 follows immediately from (2.7). □

It thus remains to prove Lemma 2.

PROOF OF LEMMA 2. It will be useful to appeal to the following setup. Consider *two* rows of $s_n + k_n$ spaces. In the first row fill the spaces from left to right with the numbered cards in increasing order. Leave the second row of spaces empty, but number its positions from left to right. Randomly choose $k_n$ spaces in the first row. The cards occupying those spaces will constitute the randomly selected $k_n$ cards in the formulation of the lemma. Now randomly choose $k_n$ spaces from the second row. Insert the $k_n$ cards chosen above into these $k_n$ spaces in the second row in increasing order. Now take the remaining $s_n$ cards from the first row and insert them into the remaining $s_n$ spaces in the second row in increasing order. Clearly, the state of the second row has the same distribution as does the state of the row described in the lemma. We will construct an increasing subsequence from this second row that uses all of the group of $s_n$ cards noted above and some from among the group of $k_n$ cards noted above.

Let $\{X_j\}_{j=1}^{k_n}$ denote the $k_n$ spaces selected from the first row and let $\{Y_j\}_{j=1}^{k_n}$ denote the $k_n$ spaces selected from the second row. In both cases we label the finite sequences so that they are increasing in $j$. Now define $Z_j = 1_{\{X_j = Y_j\}}$ and $T_n = \sum_{j=1}^{s_n + k_n} Z_j$. A moment's thought reveals that there exists an increasing subsequence of length $s_n + T_n$ from the second row. For technical reasons, we will work with $\hat{T}_n = \sum_{j=[(1/4)k_n]+1}^{[(3/4)k_n]-1} Z_j$. We will denote expectations involving $\hat{T}_n$ by $\mathcal{E}$, since $E$ has been reserved for the expectation with respect to the uniform measure on $S_n$. We will prove the following inequalities:

LEMMA 3. *Let $s_n$ and $k_n$ be as in Lemma* 2. *Then*

$$(2.8) \qquad \mathcal{E}\hat{T}_n \geq C_1 \frac{k_n^{3/2}}{n^{1/2}} \qquad \text{for some } C_1 > 0.$$

LEMMA 4. *Let $s_n$ and $k_n$ be as in Lemma 2. Then*

$$\mathcal{E}\hat{T}_n^2 \leq C_2 \frac{k_n^3}{n} \qquad \text{for some } C_2 > 0. \tag{2.9}$$

From (2.8) and (2.9), a standard argument gives

$$\liminf_{n\to\infty} \text{Prob}\left(\hat{T}_n > \delta \frac{k_n^{3/2}}{n^{1/2}}\right) > 0 \qquad \text{for some } \delta > 0. \tag{2.10}$$

Indeed, assume to the contrary that (2.10) does not hold. Then taking a subsequence if necessary and using (2.8), we may assume without loss of generality that $\frac{\hat{T}_n}{\mathcal{E}(\hat{T}_n)}$ converges to 0 in probability. However, by (2.8) and (2.9), $\mathcal{E}(\frac{\hat{T}_n}{\mathcal{E}(\hat{T}_n)})^2$ is bounded; thus, $\frac{\hat{T}_n}{\mathcal{E}(\hat{T}_n)}$ is uniformly integrable. The uniform integrability along with the convergence to 0 in probability force one to conclude that $\mathcal{E}\frac{\hat{T}_n}{\mathcal{E}(\hat{T}_n)}$ converges to 0, which is a contradiction since this expectation equals 1 for all $n$. Lemma 2 follows immediately from (2.10) and the fact that there exists an increasing subsequence of length $s_n + T_n$. □

REMARK. We suspect that the law of large numbers holds for $\hat{T}_n$; that is, $\lim_{n\to\infty} \frac{\hat{T}_n}{\mathcal{E}\hat{T}_n} = 1$ in probability. Using this, it's easy to refine the above argument to show that the law of large numbers does not hold for $Z_{n;k_n}$ when $k_n = n^{4/9}$.

It remains to prove Lemmas 3 and 4. We begin with Lemma 3.

PROOF OF LEMMA 3. In order to simplify the notation and render the proof more transparent, we will prove the lemma under the assumption that $k_n = n^l$, for some $l \in (\frac{1}{3}, \frac{1}{2})$, and $s_n + k_n = n^{1/2}$. The same method works with $s_n$ and $k_n$ as in the statement of Lemma 3, as long as $k_n = o(n^{1/2})$. If this condition on $k_n$ is not satisfied, then the structure of the proof must be amended slightly. We leave this to the reader since the important case is when $l$ is close to $\frac{4}{9}$. Note that for $1 \leq j \leq n^l$ and $j \leq r \leq n^{1/2}$, we have

$$\text{Prob}(X_j = r) = \frac{\binom{r-1}{j-1}\binom{n^{1/2}-r}{n^l-j}}{\binom{n^{1/2}}{n^l}} = \frac{j}{r}\frac{\binom{r}{j}\binom{n^{1/2}-r}{n^l-j}}{\binom{n^{1/2}}{n^l}}. \tag{2.11}$$

We calculate $P(X_j = r)$, for $j = \gamma k_n = \gamma n^l$ and $r = \gamma(s_n + k_n) + \delta k_n^{q/l} = \gamma n^{1/2} + \delta n^q$, where $\gamma \in [\frac{1}{4}, \frac{3}{4}]$, $q \in [0, \frac{1-l}{2}]$ and $\delta \in [-1, 1]$. Furthermore, given $l$ and $q$, $\gamma$ and $\delta$ are chosen so that $j$ and $r$ are integers. Opening up the binomial terms, using Stirling's approximation, doing some algebra and making cancellations, we find that

$$\frac{\binom{r}{j}\binom{n^{1/2}-r}{n^l-j}}{\binom{n^{1/2}}{n^l}}$$

$$\sim (2\pi\gamma(1-\gamma)n^l)^{-1/2}\left(1+\frac{\delta}{\gamma}n^{q-1/2}\right)^{(\gamma n^{1/2}+\delta n^q)}$$

(2.12)
$$\times \left(1-\frac{\delta}{1-\gamma}n^{q-1/2}\right)^{((1-\gamma)n^{1/2}-\delta n^q)}(1-n^{l-1/2})^{(n^{1/2}-n^l)}$$

$$\times \left(1-n^{l-1/2}-\frac{\delta}{1-\gamma}n^{q-1/2}\right)^{(-(1-\gamma)n^{1/2}+(1-\gamma)n^l+\delta n^q)}$$

$$\times \left(1+\frac{\delta}{\gamma}n^{q-1/2}-n^{l-1/2}\right)^{(-\gamma n^{1/2}-\delta n^q+\gamma n^l)} \quad \text{as } n \to \infty.$$

In order to evaluate the asymptotic behavior of the product of the final five factors on the right-hand side of (2.12), we take the logarithm of this product, obtaining

$$(\gamma n^{1/2}+\delta n^q)\log\left(1+\frac{\delta}{\gamma}n^{q-1/2}\right)$$

$$+((1-\gamma)n^{1/2}-\delta n^q)\log\left(1-\frac{\delta}{1-\gamma}n^{q-1/2}\right)$$

(2.13)
$$+(n^{1/2}-n^l)\log(1-n^{l-1/2})$$

$$+(-(1-\gamma)n^{1/2}+(1-\gamma)n^l+\delta n^q)\log\left(1-n^{l-1/2}-\frac{\delta}{1-\gamma}n^{q-1/2}\right)$$

$$+(-\gamma n^{1/2}-\delta n^q+\gamma n^l)\log\left(1+\frac{\delta}{\gamma}n^{q-1/2}-n^{l-1/2}\right).$$

Note that the logarithmic terms above, whose arguments differ from 1 by a linear combination of $n^{q-1/2}$ and $n^{l-1/2}$, are multiplied by a linear combination of positive powers of $n$. We will now use the Taylor series expansion $\log(1-x) = -x - \frac{x^2}{2} - \frac{x^3}{3} - \cdots$ to investigate the behavior of (2.13). If the highest power of $n$ remaining after cancellations in the expansion of (2.13) is non-positive, then it follows that the product of the final five factors on the right-hand side of (2.12) is on the order of unity. On the other hand, if this highest power is positive, then the coefficient of this term is necessarily negative [if it were positive, one would have $\text{Prob}(X_j = r) > 1$], and thus this product of five factors decays super-polynomially. We will show that, for the ranges of $l$ and $q$ as specified at the beginning of the proof of the lemma, the first possibility obtains. Notice that $l$ appears only in the final three terms in (2.13), and when $\delta$ is set equal to 0 (which has the effect of deleting $q$), the sum of these three terms is identically equal to 0. Thus, it follows that when one makes the expansion above, all the terms involving powers of $n$ which do not depend on $q$ will cancel out. Consequently, we need only consider powers of $n$ that depend on $q$.

We show now that it suffices to expand the logarithms in (2.13) up to order two. Note that conditions on $l$ and $q$ guarantee that $q < \frac{1}{3}$. Consider the first term on the right-hand side of (2.13). The contribution here from the third-order term in the expansion of the logarithm is on the order $n^{3(q-1/2)}$. The larger of the two terms multiplying this logarithmic term is $n^{1/2}$. Consequently, in the first term on the right-hand side of (2.13), the contribution from the third-order term in the expansion of the logarithm is on the order of $n^{1/2+3(q-1/2)}$. The exponent here is negative since $q < \frac{1}{3}$. This analysis shows that we need only expand the first term on the right-hand side of (2.13) up to the second order terms. The same argument works for the second term on the right-hand side of (2.13). The third term on the right-hand side of (2.13) can be ignored since it does not depend on $q$. Now consider the fourth and fifth terms on the right-hand side of (2.13) together. The $m$th-order terms in the expansion of the two logarithms will be

$$
(2.14) \quad \begin{aligned}
&-\frac{1}{m}\left(n^{l-1/2} + \frac{\delta}{1-\gamma}n^{q-1/2}\right)^m \\
&= -\frac{1}{m}\sum_{j=0}^{m}\binom{m}{j}n^{(m-j)(l-1/2)}\left(\frac{\delta}{1-\gamma}n^{q-1/2}\right)^j
\end{aligned}
$$

and

$$
-\frac{1}{m}\left(n^{l-1/2} - \frac{\delta}{\gamma}n^{q-1/2}\right)^m = -\frac{1}{m}\sum_{j=0}^{m}\binom{m}{j}n^{(m-j)(l-1/2)}\left(-\frac{\delta}{\gamma}n^{q-1/2}\right)^j.
$$

Recalling (2.13), one must multiply these two terms above respectively by $-(1-\gamma)n^{1/2} + (1-\gamma)n^l + \delta n^q$ and $-\gamma n^{1/2} - \delta n^q + \gamma n^l$. First consider the contribution on the right-hand side of (2.14) coming from $j = 0$. In this case, there is no dependence on $q$. Since we are only interested in powers of $n$ that depend on $q$, we need only consider multiplying these two terms respectively by $\delta n^q$ and $-\delta n^q$. Consequently, these contributions cancel each other out. Now consider the case $j = 1$. In this case, the contributions from multiplying by $n^{1/2}$ and by $n^l$ all cancel out. The contributions from multiplying by $n^q$ give a term on the order $n^{(m-1)(l-1/2)+(q-1/2)+q}$. Since $q \leq \frac{1-l}{2}$, it follows that the exponent of $n$ here is nonpositive for all $m \geq 2$. Now consider the case $j = 2$. Here for the first time there are no cancellations.

Upon multiplying things out, the highest-order term will be on the order of $n^{(m-2)(l-1/2)+2(q-1/2)+1/2}$. Using the fact that $q \le \frac{1-l}{2}$, it follows that the exponent of $n$ here is nonpositive for $m \ge 3$. Since the conditions on $l$ and $q$ force $q < l$, it follows a fortiori that all powers of $n$ appearing after multiplication will be nonpositive whenever $j \ge 2$ and $m \ge 3$. From the above analysis, we conclude that it suffices to expand all the logarithms on the right side of (2.13) up to order 2 and ignore all the higher-order terms.

If one considers the contribution from the first-order terms [i.e., replacing $\log(1-x)$ with $-x$ in (2.13)], one finds that everything cancels out. And when one considers the contribution in (2.13) from the second-order terms [i.e., replacing $\log(1-x)$ with $-\frac{x^2}{2}$], one finds that all but one of the terms cancel, leaving $-\frac{3}{2}\frac{\delta^2}{\gamma(1-\gamma)}n^{2q+l-1}$. Since we have assumed that $q \le \frac{1-l}{2}$, it follows that $2q+l-1 \le 0$. We have thus shown that the product of the final five terms on the right-hand side of (2.12) is on the order of unity as long as $q \le \frac{1-l}{2}$. We conclude then from (2.11) and (2.12) and the fact that $\gamma \in [\frac{1}{4}, \frac{3}{4}]$ that there exists a $c_0 > 0$ such that $\text{Prob}(X_j = r) \ge c_0 n^{(l-1)/2}$. Since this inequality is in force for every $q \le \frac{1-l}{2}$ and every $\delta \in [-1, 1]$, it follows that, for each $\gamma$ (and consequently, each $j$), there are $[2n^{(1-l)/2}]$ different values of $r$ for which this inequality holds. Thus $\text{Prob}(Z_j = 1) = \text{Prob}(X_j = Y_j) \ge [2n^{(1-l)/2}](c_0 n^{(l-1)/2})^2 = c_1 n^{(l-1)/2}$, for some $c_1 > 0$. As $\gamma$ varies from $\frac{1}{4}$ to $\frac{3}{4}$, we obtain at least $[\frac{1}{2}n^l] - 2$ different values of $j$ between $[\frac{1}{4}k_n]+1$ and $[\frac{3}{4}k_n]-1$ for which the above inequality holds. Thus, $\mathcal{E}\hat{T}_n = \sum_{j=[(1/4)k_n]+1}^{[(3/4)k_n]-1} \mathcal{E}Z_j = \sum_{j=[(1/4)k_n]+1}^{[(3/4)k_n]-1} \text{Prob}(Z_j=1) > \frac{1}{3}n^l(c_1 n^{(l-1)/2}) = \frac{1}{3}c_1 n^{(3l-1)/2} = \frac{1}{3}c_1 \frac{k_n^{3/2}}{n^{1/2}}$. □

PROOF OF LEMMA 4. In the proof of Lemma 3 we assumed that $k_n = n^l$ and $s_n + k_n = n^{1/2}$. For the proof here, we return to the original notation. We have

$$\mathcal{E}\hat{T}_n^2 = 2 \sum_{[(1/4)k_n]+1 \le i < j \le [(3/4)k_n]-1} \text{Prob}(Z_i = Z_j = 1)$$

(2.15)

$$+ \sum_{j=[(1/4)k_n]+1}^{[(3/4)k_n]-1} \text{Prob}(Z_j = 1).$$

We will use the following bound, whose proof we postpone.

(2.16) $\text{Prob}(X_j = r) \le \begin{cases} C\frac{k_n}{j^{1/2}(s_n+k_n)}, & \text{for } 1 \le j \le \frac{k_n}{2} \text{ and all } r, \\ C\frac{k_n}{(k_n-j+1)^{1/2}(s_n+k_n)}, & \text{for } \frac{k_n}{2} < j \le k_n \text{ and all } r. \end{cases}$

We have

(2.17) $\text{Prob}(Z_j = 1) = \sum_{r=1}^{s_n+k_n} \text{Prob}(X_j = Y_j = r) = \sum_{r=1}^{s_n+k_n} (\text{Prob}(X_j = r))^2.$

In order to use (2.16) to obtain an upper bound for (2.17), consider the maximum value attained by the expression $\sum_{j=1}^{m} p_j^2$, where $\{p_j\}_{j=1}^{m}$ is subject to the constraints $\sum_{j=1}^{m} p_j = 1$ and $0 \leq p_j \leq \delta$, for all $j$. (Of course, it necessarily follows that $\delta \geq \frac{1}{m}$.) The maximum is achieved by letting $p_j = \delta$, for $j = 1, \ldots, m_0$, where $m_0$ is defined by $m_0 \delta \leq 1$ and $(m_0 + 1)\delta > 1$, and then letting $p_{m_0+1} = 1 - m_0 \delta$ and $p_j = 0$ for all $j > m_0 + 1$. It is easy to check that the maximum value is no greater than $\delta$. Applying this to (2.16) and (2.17), it follows that

$$(2.18) \quad \operatorname{Prob}(Z_j = 1) \leq \begin{cases} C \frac{k_n}{j^{1/2}(s_n + k_n)}, & \text{if } 1 \leq j \leq \frac{k_n}{2}, \\ C \frac{k_n}{(k_n - j + 1)^{1/2}(s_n + k_n)}, & \text{if } \frac{k_n}{2} < j \leq k_n. \end{cases}$$

Thus,

$$(2.19) \quad \sum_{j=1}^{k_n} \operatorname{Prob}(Z_j = 1) \leq C_1 \frac{k_n^{3/2}}{s_n + k_n}.$$

In light of (2.15) and (2.19), to complete the proof of the lemma it remains to show that

$$(2.20) \quad \sum_{[(1/4)k_n]+1 \leq i < j \leq [(3/4)k_n]-1} \operatorname{Prob}(Z_i = Z_j = 1) = O\left(\frac{k_n^3}{(s_n + k_n)^2}\right).$$

Recalling that $l > \frac{1}{3}$, let $\varepsilon > 0$ be sufficiently small so that $\frac{1-l}{2} + \varepsilon < \frac{1}{2}$. For $i < j$, we write

$$\operatorname{Prob}(Z_i = Z_j = 1)$$

$$= \operatorname{Prob}\left(Z_i = 1, X_i \leq \frac{i}{k_n}(s_n + k_n) + n^{(1-l)/2+\varepsilon}\right)$$

$$(2.21)$$

$$\times \operatorname{Prob}\left(Z_j = 1 | Z_i = 1, X_i \leq \frac{i}{k_n}(s_n + k_n) + n^{(1-l)/2+\varepsilon}\right)$$

$$+ \operatorname{Prob}\left(Z_i = 1, Z_j = 1, X_i > \frac{i}{k_n}(s_n + k_n) + n^{(1-l)/2+\varepsilon}\right).$$

The proof of Lemma 3, in particular the last paragraph where the condition $q \leq \frac{1-l}{2}$ was used, shows that $\operatorname{Prob}(X_i > \frac{i}{k_n}(s_n + k_n) + n^{(1-l)/2+\varepsilon})$ decays to zero super-polynomially as $n \to \infty$. (Here $\frac{i}{k_n}$ fulfills the role of $\gamma$ in the proof of Lemma 3.) In light of this fast decay, it follows from (2.21) that (2.20) will hold if we show that

$$\sum_{[(1/4)k_n]+1 \leq i < j \leq [(3/4)k_n]-1} \operatorname{Prob}(Z_i = 1)$$

$$(2.22) \qquad \times P\left(Z_j = 1 | Z_i = 1, X_i \leq \frac{i}{k_n}(s_n + k_n) + n^{(1-l)/2+\varepsilon}\right)$$

$$= O\left(\frac{k_n^3}{(s_n + k_n)^2}\right).$$

Using the structure of $\{Z_i\}_{i=1}^{k_n}$, we can apply (2.18) to bound from above $\text{Prob}(Z_j = 1 | Z_i = 1, X_i \leq \frac{i}{k_n}(s_n + k_n) + n^{(1-l)/2+\varepsilon})$. We must replace $k_n$, $j$ and $(s_n + k_n)$, respectively, in (2.18) by $k_n - i$, $j - i$ and $(s_n + k_n - \frac{i}{k_n}(s_n + k_n) - n^{(1-l)/2+\varepsilon})$. Thus,

(2.23)
$$\text{Prob}\left(Z_j = 1 | Z_i = 1, X_i \leq \frac{i}{k_n}(s_n + k_n) + n^{(1-l)/2+\varepsilon}\right)$$
$$\leq \begin{cases} C\frac{k_n - i}{(j-i)^{1/2}(s_n + k_n - (i/k_n)(s_n + k_n) - n^{(1-l)/2+\varepsilon})}, \\ \quad \text{if } 1 \leq j - i \leq \frac{k_n - i}{2}, \\ C\frac{k_n - i}{(k_n - j + 1)^{1/2}(s_n + k_n - (i/k_n)(s_n + k_n) - n^{(1-l)/2+\varepsilon})}, \\ \quad \text{if } \frac{k_n - i}{2} < j - i \leq k_n - i. \end{cases}$$

Since we have chosen $\varepsilon$ such that $\frac{1-l}{2} + \varepsilon < \frac{1}{2}$ and since $s_n + k_n$ is on the order of $n^{1/2}$, it follows from (2.23) that

(2.24)
$$\text{Prob}\left(Z_j = 1 | Z_i = 1, X_i \leq \frac{i}{k_n}(s_n + k_n) + n^{(1-l)/2+\varepsilon}\right)$$
$$\leq \begin{cases} C_1 \frac{k_n}{(j-i)^{1/2}(s_n + k_n)}, & \text{if } i+1 \leq j \leq \frac{k_n+i}{2} \text{ and } \frac{1}{4}k_n \leq i, j \leq \frac{3}{4}k_n, \\ C_1 \frac{k_n^{1/2}}{s_n + k_n}, & \text{if } j > \frac{k_n+i}{2} \text{ and } \frac{1}{4}k_n \leq i, j \leq \frac{3}{4}k_n. \end{cases}$$

From (2.18) and (2.24), we obtain

(2.25)
$$\text{Prob}(Z_i = 1)\text{Prob}\left(Z_j = 1 | Z_i = 1, X_i \leq \frac{i}{k_n}(s_n + k_n) + n^{(1-l)/2+\varepsilon}\right)$$
$$\leq \begin{cases} C_2 \frac{k_n^{(3/2)}}{(j-i)^{(1/2)}(s_n + k_n)^2}, \\ \quad \text{if } i+1 \leq j \leq \frac{k_n+i}{2} \text{ and } \frac{1}{4}k_n \leq i, j \leq \frac{3}{4}k_n, \\ C_2 \frac{k_n}{(s_n + k_n)^2}, \\ \quad \text{if } j > \frac{k_n+i}{2} \text{ and } \frac{1}{4}k_n \leq i, j \leq \frac{3}{4}k_n. \end{cases}$$

Using (2.25) we obtain

(2.26)
$$\text{Prob}(Z_i = 1) \sum_{j=i+1}^{[(3/4)k_n]-1} \text{Prob}\left(Z_j = 1 | Z_i = 1, X_i \leq \frac{i}{k_n}(s_n + k_n) + n^{(1-l)/2+\varepsilon}\right)$$
$$\leq C_3 \frac{k_n^2}{(s_n + k_n)^2} \quad \text{for } \frac{1}{4}k_n \leq i \leq \frac{3}{4}k_n.$$

Summing (2.26) over $i$ gives (2.20), which completes the proof of the lemma.

We now return to prove (2.16). By symmetry, the second inequality in (2.16) follows from the first one. Thus, we consider only the first one. In the present notation, (2.11) becomes

$$\text{Prob}(X_j = r) = \frac{\binom{r-1}{j-1}\binom{s_n+k_n-r}{k_n-j}}{\binom{s_n+k_n}{k_n}} = \frac{j}{r}\frac{\binom{r}{j}\binom{s_n+k_n-r}{k_n-j}}{\binom{s_n+k_n}{k_n}}. \tag{2.27}$$

One can check readily that $\text{Prob}(X_1 = r)$ achieves its maximum value at $r = 1$, where it equals $\frac{k_n}{s_n+k_n}$. Thus, from now on, we will consider $\text{Prob}(X_j = r)$ under the assumption that $j \geq 2$. Consider $H(r) \equiv \binom{r-1}{j-1}\binom{s_n+k_n-r}{k_n-j}$. Writing out $\frac{H(r)}{H(r-1)}$ explicitly and solving for $r$ (considered now as a real number) in the equation $\frac{H(r)}{H(r-1)} = 1$, one finds that $r = 1 + \frac{(j-1)(s_n+k_n)}{k_n-1}$. This shows that the function $H(r)$ obtains its maximum at an integer in the interval $[\frac{(j-1)(s_n+k_n)}{k_n-1}, \frac{(j-1)(s_n+k_n)}{k_n-1} + 2]$. Denote an integer where the maximum occurs by $r_0$, suppressing the dependence on $n$ and $j$. By Stirling's approximation, there exist constants $C_1, C_2 > 0$ such that $C_1 m^m \exp(-m) m^{1/2} \leq m! \leq C_2 m^m \exp(-m) m^{1/2}$, for all $m \geq 1$. Using this two-sided bound, we have

$$\frac{\binom{r_0}{j}\binom{s_n+k_n-r_0}{k_n-j}}{\binom{s_n+k_n}{k_n}}$$

$$\leq C\left(\frac{r_0(s_n+k_n-r_0)k_n s_n}{j(r_0-j)(k_n-j)(s_n-r_0+j)(s_n+k_n)}\right)^{1/2} \tag{2.28}$$

$$\times \frac{r_0^{r_0}(s_n+k_n-r_0)^{(s_n+k_n-r_0)}k_n^{k_n}s_n^{s_n}}{j^j(r_0-j)^{(r_0-j)}(k_n-j)^{(k_n-j)}(s_n-r_0+j)^{(s_n-r_0+j)}(s_n+k_n)^{(s_n+k_n)}},$$

for some $C > 0$. The conditions that have been placed on $s_n, k_n, r_0$ and $j$ guarantee that

$$\frac{r_0}{r_0-j}\frac{k_n}{k_n-j}\frac{s_n+k_n-r_0}{s_n-r_0+j}\frac{s_n}{s_n+k_n} \quad \text{is bounded from above.} \tag{2.29}$$

Furthermore, we claim that

$$\frac{r_0^{r_0}(s_n+k_n-r_0)^{(s_n+k_n-r_0)}k_n^{k_n}s_n^{s_n}}{j^j(r_0-j)^{(r_0-j)}(k_n-j)^{(k_n-j)}(s_n-r_0+j)^{(s_n-r_0+j)}(s_n+k_n)^{(s_n+k_n)}} \leq 1. \tag{2.30}$$

This follows from the following lemma.

LEMMA 5. *For $a, b, c, d > 0$, one has*

$$\frac{(a+b)^{(a+b)}}{a^a b^b}\frac{(c+d)^{(c+d)}}{c^c d^d}\frac{(a+c)^{(a+c)}(b+d)^{(b+d)}}{(a+b+c+d)^{(a+b+c+d)}} \leq 1. \tag{2.31}$$

We postpone the proof of Lemma 5 until the completion of the proof of Lemma 4. The proof of (2.30) follows from Lemma 5 by choosing $a = j, b = r_0 - j, c = k_n - j$ and $d = s_n - r_0 + j$. Now (2.16) follows from (2.27), (2.28), (2.29), (2.30) and the definition of $r_0$. This completes the proof of Lemma 4. □

It remains to prove Lemma 5.

PROOF OF LEMMA 5. We will show that the logarithm of the left-hand side of (2.31) is nonpositive. The logarithm of the left-hand side of (2.31) can be written as $F(a,b) + F(c,d) - F(a+c, b+d)$, where

$$F(x,y) = (x+y)\log(x+y) - x\log x - y\log y.$$

Thus, we need to show that

(2.32) $$F(a+c, b+d) - F(a,b) - F(c,d) \geq 0.$$

To prove (2.32), it suffices to show that $G(t) \equiv F(a+c, b+t) - F(a,b) - F(c,t)$ is nonnegative for all $t \geq 0$. Differentiating, one finds that $G$ has only one critical point and that $G$ attains its minimum there. This critical point is at $t = \frac{bc}{a}$. Thus, it only remains to show that $G(\frac{bc}{a}) \geq 0$. A direct calculation reveals that $G(\frac{bc}{a}) = 0$. □

DEPARTMENT OF MATHEMATICS
TECHNION—ISRAEL INSTITUTE OF TECHNOLOGY
HAIFA 32000
ISRAEL
E-MAIL: pinsky@math.technion.ac.il
URL: http://www.math.technion.ac.il/~pinsky/